\documentclass[a4paper,11pt,reqno]{amsart}
\usepackage[latin1]{inputenc}

\usepackage{amsmath}
\usepackage{graphicx,color}
\usepackage{bbm}
\usepackage{xspace}

\definecolor{darkblue}{rgb}{.15,.15,.8} 
\usepackage[colorlinks,
            linkcolor=darkblue, 
            citecolor=darkblue,
            pagecolor=darkblue,
            filecolor=darkblue,
            urlcolor=darkblue,
            pdfdisplaydoctitle,
            anchorcolor=darkblue,
            menucolor=darkblue,
            dvips]{hyperref}
\usepackage[textwidth=.72\paperwidth]{geometry}
\usepackage{amsfonts, amsmath, amssymb, amsthm}
\usepackage[mathscr]{eucal}

\hypersetup{
pdfauthor={Michael Joswig, Bejamin M\"uller, Andreas Paffenholz},
pdftitle={polymake and Lattice Polytopes},
pdfkeywords={polymake system, lattice polytope, Hilbert basis, toric geometry}}

\providecommand{\ZZ}{{\mathbb Z}}
\providecommand{\RR}{{\mathbb R}}
\providecommand{\cI}{{\mathscr{I}}}
\providecommand{\cN}{{\mathscr{N}}}

\providecommand\polymake{\texttt{polymake}\xspace}
\providecommand\normaliz{\texttt{normaliz2}\xspace}
\providecommand\latte{\texttt{Latte macchiato}\xspace}
\providecommand\fourtytwo{\texttt{4ti2}\xspace}
\providecommand\nauty{\texttt{nauty}\xspace}

\newenvironment{dense_itemize}{%
\begin{list}{$\triangleright$}%
{\setlength{\topsep}{1mm}%
\setlength{\partopsep}{0mm}%
\setlength{\parskip}{0mm}%
\setlength{\parsep}{0mm}%
\setlength{\itemsep}{0mm}%
\setlength{\labelwidth}{4mm}%
\setlength{\leftmargin}{0mm}%
\addtolength{\leftmargin}{\labelwidth}%
\addtolength{\leftmargin}{\labelsep}%
\setlength{\itemindent}{0mm}}}%
{\end{list}}

\newtheorem{theorem}{Theorem}

\title{polymake and Lattice Polytopes}

\markleft{Joswig, M\"uller\& Paffenholz}

\author{Michael Joswig}
\address{Institut f\"ur Mathematik, MA 6-2, TU Berlin, 10623 Berlin, Germany}
\email{joswig@math.tu-berlin.de}
\thanks{Research of Michael Joswig supported by German Research
  Foundation (DFG) Research Unit ``Polyhedral Surfaces''.}
\author{Benjamin M\"uller}
\address{FU Berlin, Institut f\"ur  Mathematik, Arnimallee 3, 14195 Berlin, Germany}
\email{benmuell@math.fu-berlin.de}
\thanks{Benjamin M\"uller is supported by Christian Haase's Emmy
  Noether fellowship HA 4383/1 of the DFG}
\author{Andreas Paffenholz}
\address{FU Berlin, Institut f\"ur  Mathematik, Arnimallee 3, 14195  Berlin, Germany}
\email{paffenho@math.fu-berlin.de}

\keywords{polymake system, lattice polytope, Hilbert basis, toric
  geometry}
\subjclass[2000]{52B20, 52-04, 14M25}

\date{\today}

\begin{document}

\begin{abstract}
  The \polymake software system deals with convex polytopes and related objects from
  geometric combinatorics.  This note reports on a new implementation of a subclass for
  lattice polytopes.  The features displayed are enabled by recent changes to the
  \polymake core, which will be discussed briefly.
\end{abstract}

\maketitle

\section{Introduction}

\polymake is a software system designed for analyzing convex polytopes, finite simplicial
complexes, graphs, and other objects.  While the system exists for more than a
decade~\cite{DMV:polymake} it was continuously developed and expanded.  The most recent
version fundamentally changes the way to interact with the system.  It now offers an
interface which looks similar to many computer algebra systems.  However, on the technical
level \polymake differs from most mathematical software systems: rule based computations
and an extendible dual Perl/C++ interface are the most important characteristics.

\polymake can now also handle hierarchies of objects where each level may come with
additional sets of rules.  \polymake handles casts between subclasses based on property
requests.  We will explain this feature by means of a new subclass
\texttt{LatticePolytope} that is derived from the existing class
\texttt{Polytope<Rational>}.  However, some of the new functions can also be applied to
any rational polytope.  A \emph{lattice polytope} is a polytope whose vertices are
contained in a lattice $\Lambda\subset\RR^n$ \cite{BeckRobbins,Barvinok08}.  \polymake
always assumes $\Lambda=\ZZ^n$.  This subclass reflects a new use of \polymake in toric
geometry, where lattice polytopes encode properties of toric varieties and toric
ideals~\cite{Oda88,Ewald}.  String theorists have been interested in special lattice
polytopes, as they led to the construction of mirror pairs of Calabi-Yau
varieties~\cite{Batyrev94}.  Gr\"obner bases of toric ideals have been applied to
optimization problems~\cite{Sturmfels96}.

We will explain all relevant concepts for our exposition on the way.  Lattice polytopes
have also become an important subject in other areas of mathematics.  Enumerating
non-negative solutions of Diophantine equations can be interpreted as counting lattice
points in a polyhedron~\cite{Stanley96}.  Contingency tables in statistics can be modeled
by lattice polytopes~\cite{DS98}.  Sampling then corresponds to finding integral points in the
polytope.

The paper is organized as follows.  First we will review the recent changes and the new
\polymake interface. Then we will report on our new implementation of a subclass for
lattice polytopes.  In particular, this comprises interfaces to \fourtytwo \cite{4ti2},
\latte \cite{latte,lattemacchiato}, and \normaliz \cite{normaliz2}.  We will show how the
user can benefit from the common interface to these systems via \polymake and how one can
extend their functionality by combining with \polymake's features.  Rather than discussing
implementation details we will explain the functions available with one easy running
example.  This note then concludes with a final section analyzing a specific
$6$-dimensional polyhedral cone which was found to be a counter-example to a conjecture of
Seb\H{o} \cite{Sebo} by Bruns et al.~\cite{BGHW99}.

\section{polymake -- the Next Generation}

The general ideas which lead to the design and the implementation of the \polymake system
more than ten years ago are still valid.  The key goals are the following.
\begin{dense_itemize}
\item The system should be scalable with the user's ability to write programs.  This means
  that basic usage should not require any programming skills, while it should be powerful
  enough not to restrain the programming expert.
\item The system should not try to ``re-invent the wheel''. There is a multitude of
  valuable pieces of software for individual tasks; so they should be suitably interfaced
  rather than their functionality be duplicated.
\item The system should be really easy to extend.  It should be possible to model new
  mathematical objects and to integrate them into the existing framework.
\end{dense_itemize}
These ``golden rules'' are most natural, and most users of mathematical software systems
would probably agree that all of these are very desirable.  For instance, the SAGE system
is following a similar strategy albeit on a somewhat larger scale \cite{sage}.  In
\polymake we are focusing on convex polytopes and related objects from the realm of
geometric combinatorics.  The ``golden rules'' already have a number of implications, some
obvious and some less obvious.  The most important design decisions which can be derived
are: The system requires both a compiled and an interpreted programming language (we
settled for C++ and Perl), and the system must be an Open Source project (we settled for
the GNU Public License).  By far the most difficult to accomplish is the third rule.  And,
in fact, a large part of \polymake's code evolution over the last decade can be seen as an
attempt to re-interpret this rule again and again with an increasing level of abstraction.

A word of warning to the experienced \polymake user.  On a technical level the new version
of \polymake is very different from previous versions.  From the point of view of the
working mathematician this results in a number of benefits.  In particular, the overall
usability is improved, while we gained additional flexibility and speed.  The unavoidable
drawback is that the interface had to be changed in a substantial way.

Using \polymake now means to start a program named ``\polymake'' from the command line,
and then to work in a shell-type environment typical for most computer algebra systems.
The language for interacting with the system is Perl, but we added a few features in order
to easy the usability.  We give a very brief overview of how to get started with the new
system.

\begin{footnotesize}
\begin{verbatim}
Welcome to polymake version 2.9.6, rev. 9033
Copyright (c) 1997-2009
Ewgenij Gawrilow (TU Berlin), Michael Joswig (TU Darmstadt)
http://www.math.tu-berlin.de/polymake,  mailto:polymake@math.tu-berlin.de

This is free software licensed under GPL; see the source for copying conditions.
There is NO warranty; not even for MERCHANTABILITY or FITNESS FOR A PARTICULAR PURPOSE.

Type 'help;' for basic instructions.

Application polytope uses following third-party software (for details: help 'credits';)
4ti2, azove, cddlib, lrslib, nauty, normaliz2, porta, qhull, splitstree, topcom, vinci
polytope > 
\end{verbatim}
\end{footnotesize}

By now there are several different \emph{applications} known to \polymake.  Each
application comprises a main \emph{object} type, \emph{properties} which describe an
object of this type, and a set of \emph{rules}.  By default the first application to start
is the one dealing with convex polytopes, and this is made visible by showing the command
line prompt ``\texttt{polytope >}''.  The last line before the prompt lists the programs
whose interfaces are loaded.  Since everything (the application, the objects, the
properties, the rules, the interfaces, and the defaults) can be modified or extended by
the user, what shows up exactly very much depends on the local installation.  The main
purpose of this note is to explain how a new sub-type for lattice polytopes is organized
within the object hierarchy for general polytopes.

The following simple example can explain the \polymake concept in a nutshell.  The first
command produces a $3$-dimensional cube with $\pm1$-coordinates (and assigns it to the
variable \verb+$P+), while the second one (separated by ``\texttt{;}'') prints its
$f$-vector, that is, the number of faces per dimension.  Clearly we have eight vertices,
12 edges, and six facets.

\begin{footnotesize}
\begin{verbatim}
polytope > $P=cube(3); print $P->F_VECTOR;
8 12 6
\end{verbatim}
\end{footnotesize}

The function \texttt{cube} returns a polytope object of type \texttt{Polytope<Rational>},
and \verb+F_VECTOR+ is a property of this class, which models polytopes with rational
coordinates.  Notice that there are polytopes whose combinatorial type does not admit any
rational representation \cite[\S6.5]{Ziegler}.  \polymake reduces computing the $f$-vector
of this cube to finding a suitable sequence of rules and to execute them one after
another.  These rules can be shown as follows.  To this end we restart from scratch.

\begin{footnotesize}
\begin{verbatim}
polytope > $P=cube(3); print join(", ", $P->list_properties);
AMBIENT_DIM, DIM, FACETS, VERTICES_IN_FACETS, BOUNDED
polytope > print $P->type->full_name;
Polytope<Rational>
\end{verbatim}
\end{footnotesize}

Our cube \texttt{\$P} is ``born'' as an object with the five initial properties \verb+AMBIENT_DIM+, \verb+DIM+,
\verb+FACETS+, \verb+VERTICES_IN_FACETS+, \verb+BOUNDED+, all of which are redundant
except for \verb+FACETS+, which gives a description of the cube as the intersection of six
affine halfspaces.

\begin{footnotesize}
\begin{verbatim}
polytope > print $P->FACETS;
1 1 0 0
1 0 0 1
1 0 1 0
1 0 -1 0
1 -1 0 0
1 0 0 -1
\end{verbatim}
\end{footnotesize} 

Each line is a vector $(\alpha_0,\alpha_1,\dots,\alpha_d)$ representing the linear
inequality $\alpha_0+\alpha_1x_1+\dots+\alpha_dx_d\ge 0$.  The property
\verb+AMBIENT_DIM+, for instance, is the dimension of the space where our polytope lives
in, that is the number $d$, which can easily be derived from each facet by counting the
number of columns.  Now, asking for the $f$-vector means that it has to be computed from
the data given somehow.  We can look at the schedule of rules necessary to accomplish this
task.

\begin{footnotesize}
\begin{verbatim}
polytope > $schedule=$P->get_schedule("F_VECTOR");
polytope > print join("\n", $schedule->list);
HASSE_DIAGRAM : VERTICES_IN_FACETS
F_VECTOR, F2_VECTOR : HASSE_DIAGRAM
\end{verbatim}
\end{footnotesize}

Each line is one rule. Each rule has its \emph{targets} to the left of the ``\verb+:+''
and its \emph{sources} to the right.  The first line says: ``I can produce the Hasse
diagram (of the face lattice) if I know which vertex is incident with which facet''.  This
is clear since it follows from the facet that the face lattice is co-atomic, that is, each
face is the intersection of facets \cite[\S2.2
]{Ziegler}.  The second line says: ``I know
how to compute the $f$-vector (and something else that we do not care to discuss now) from
the Hasse diagram.  Each rule comes with a piece of (Perl) code which actually implements
what the rule heads shown promise.  The schedule is an object of its own right, and it can
be applied to the cube, which means that the corresponding Perl code is executed in the
order of the schedule.

\begin{footnotesize}
\begin{verbatim}
polytope > $schedule->apply($P);
polytope > print join(", ", $P->list_properties);
AMBIENT_DIM, DIM, FACETS, VERTICES_IN_FACETS, BOUNDED, HASSE_DIAGRAM, F_VECTOR,
F2_VECTOR
\end{verbatim}
\end{footnotesize}

We see that the list of properties known about our cube changed.
Three new properties have been added, and these correspond to the
total of three targets of the two rules above.  If we now ask for the
$f$-vector this information is already stored with the cube object,
and it is read from memory rather than re-computed.

As far as technology is concerned, the function \verb+cube+ which was
called to produce the cube is written in C++.  On top of the standard
Perl-interface to C we built a shared memory mechanism to access an
object from the C++ and the Perl side.  This is also fully extendible,
which means that the user is welcome to add new functions to produce
other polytopes, new properties of polytopes, or new rules to compute
existing properties in a different way.  The integration of new
functions, properties, and rules is seamless, that is, they cannot be
distinguished from the built-in ones.

There are many more things to be said about this concept both from the
logical and the technical point of view, but for the details we refer
the reader to \cite{DMV:polymake} and to further documentation at
\url{http://www.opt.tu-darmstadt.de/polymake}.

\section{Lattice Polytopes as a Subclass}

The new version of \polymake can now handle derived classes of objects specified by some
\emph{preconditions} that inherit all properties and rules from their base class but may
provide additional rules that are specific for their class.  A user may, but doesn't have
to, specify, that the object he defines falls in this class.  \polymake decides upon what
properties a user asks for, whether the object should be \emph{cast} into this subclass.
Of course, before performing the cast, \polymake checks whether the object meets the
requirements for the subclass.  The first occurrence of this new mechanism is in the class
\texttt{LatticePolytope} derived from \texttt{Polytope<Rational>}.  In \polymake, a
lattice polytope is a bounded rational polytope whose vertices are in the integer lattice
$\ZZ^n$.  The new rules in this object class concern properties of such polytopes in
connection with toric algebra and algebraic geometry.

The main focus of our implementation concerning lattice polytopes is toric geometry, so we
explain this connection here.  Let $P$ be a lattice polytope.  The \emph{normal fan}
$\cN_P$ defines a projective \emph{toric variety} $X_P$ \cite{Ewald,Sturmfels96}.  The
defining ideal $\cI_P$ of $X_p$ is a homogeneous \emph{toric ideal}.  Many properties of
the variety are reflected in the corresponding polytope.  We will see some entries in this
``dictionary'' which translates back and forth below.

There are several software packages available which proved to be useful in applications in
this area.  \normaliz by Bruns and Ichim~\cite{normaliz2} computes Hilbert bases and
$h^*$-polynomials. \latte by K\"oppe \cite{lattemacchiato} builds on previous work by De
Loera et~al.~\cite{latte}, and its key application is to count lattice points and to
compute Ehrhart polynomials.  \fourtytwo by Hemmecke et.\ al.~\cite{4ti2} solves integral
equations over $\ZZ$ and it computes convex hulls as well as Hilbert bases.  \polymake now
provides a unified access to these programs.  Additionally, we implemented various rules
to compute further important properties of lattice polytopes which can be derived.  We
will browse through the main features by using the $3$-cube from above as our running
example.  As already mentioned, our cube is the convex hull of all $\pm1$-vectors, so the
vertices do lie in the $\ZZ^3$-lattice.  We can let \polymake check this for us.

\begin{footnotesize}
\begin{verbatim}
polytope > print $P->LATTICE;
1
\end{verbatim}
\end{footnotesize}

Here the output ``$1$'' represents the boolean value ``true''.  For instance, we can ask
for the number of lattice points contained in the cube, that is, for the number
$|[-1,1]^3\cap \ZZ^3|$.  In our case, we should obtain ``$27$'' as the answer, there is
exactly one lattice point contained in the relative interior of each non-empty face.

\begin{footnotesize}
\begin{verbatim}
polytope > print $P->N_LATTICE_POINTS;
polymake: used package latte
 LattE macchiato is an improved version of LattE, a free software dedicated
 to the problems of counting and detecting lattice points inside convex polytopes,
 and the solution of integer programs.
 Copyright by Matthias Koeppe, Jesus A. De Loera and others.
 http://www.math.ucdavis.edu/~mkoeppe/latte/

27
\end{verbatim}
\end{footnotesize}

As shown \latte was called for the computation.  It uses an enhanced version of Barvinok's
algorithm\cite{koeppeirrational,barvinok_algo}.  By default \polymake gives credit to a
program when it calls it for the first time.  The corresponding output is omitted in some
of the computations below; but we will explain which package was called in each case.

\begin{footnotesize}
\begin{verbatim}
polytope > print $P->N_INTERIOR_LATTICE_POINTS;
1
\end{verbatim}
\end{footnotesize}

Sometimes it is important to know how many of the lattice points in a polytope are
contained in the interior.  While the Barvinok algorithm avoids to enumerate the points,
the user can force the complete enumeration.  This will be computed by \fourtytwo.

\begin{footnotesize}
\begin{verbatim}
print $P->INTERIOR_LATTICE_POINTS;
1 0 0 0
\end{verbatim}
\end{footnotesize}

Up to this point, none of the rules used was specific to lattice polytopes.  To the
contrary, all this makes perfect sense for any rational polytope.

Let us now switch to some properties that are only defined for lattice polytopes.  A
lattice polytope is \emph{reflexive}, if the origin is in the interior of the polytope and
all facets have integral distance one from the origin.  Equivalently, a lattice polytope
is reflexive, if also its polar is a lattice polytope.  In algebraic geometry, these
polytopes correspond to Gorenstein toric Fano varieties.  These polytopes were introduced
by Batyrev~\cite{Batyrev94} to construct mirror pairs of \emph{Calabi-Yau} varieties in
the context of string theory.  A necessary condition for a polytope to be reflexive is,
that the origin is the unique interior lattice point, so the cube is a candidate.

\begin{footnotesize}
\begin{verbatim}
polytope > print $P->REFLEXIVE;
1
\end{verbatim}
\end{footnotesize}

Of course, this is not a surprise, as the polar dual of the cube (with $\pm1$-coordinates)
is the regular octahedron, the convex hull of the standard basis vectors and their
negatives.  Reflexivity is a property that is only defined for lattice polytopes, and so
at this point, \polymake has internally cast the cube to the subclass
\texttt{LatticePolytope}.

\begin{footnotesize}
\begin{verbatim}
polytope > print $P->type->full_name;
LatticePolytope
\end{verbatim}
\end{footnotesize}

Notice the difference to the first call of the same command at the very beginning.  Many
further properties of Fano varieties can be checked via the corresponding polytope.
Reflexive polyhedra have been classified in dimensions up to $4$, see~\cite{KSA}.  In
dimension $3$, there are $124$ of these, of which $18$ correspond to smooth Fano
varieties.  The toric variety $X_P$ is \emph{smooth} (or \emph{non-singular}) if every
cone in the normal fan $\cN_P$ is unimodular.  A cone is \emph{unimodular}, if its minimal
integral generators can be extended to a basis of $\ZZ^n$.

\begin{footnotesize}
\begin{verbatim}
polytope > print $P->SMOOTH;
1
\end{verbatim}
\end{footnotesize}

We will now explore a different aspect of lattice polytopes.  Stanley
showed~\cite{Stanleyhstar} that for any $d$-dimensional polytope $P$ there is a polynomial
$h^*\in\ZZ[t]$ of degree at most $d$ with non-negative coefficients such that
\begin{align*}
  \sum_{k\ge 0}|kP\cap \ZZ^d|t^k \ = \ \frac{h^*(t)}{(1-t)^{d+1}} \, .
\end{align*}
The polynomial $h^*(t)=\sum_{k=0}^dh^*_kt^k$ is the \emph{$h^*$-polynomial} of $P$.  It is
closely related to the Ehrhart polynomial.  Some of the coefficients have a combinatorial
meaning.  For instance, $h_d^*$ counts the number of interior lattice points, while the
sum of all coefficients is the normalized volume of the polytope.  The \emph{normalized
  volume} of a $d$-dimensional polytope is $d!$ times the $d$-dimensional Euclidean
volume.  We can compute the coefficients for the cube (starting with the constant
coefficient).

\begin{footnotesize}
\begin{verbatim}
polytope > print $P->H_STAR_VECTOR;
1 23 23 1
polytope > print $P->LATTICE_VOLUME;
48
\end{verbatim}
\end{footnotesize}

\polymake calls \normaliz to compute this.  The \emph{degree} $\delta$ of $P$ is defined
as the degree of the $h^*$-polynomial.

\begin{footnotesize}
\begin{verbatim}
polytope > print $P->LATTICE_DEGREE;
3
\end{verbatim}
\end{footnotesize}

The value $d+1-\delta$ is the smallest factor by which we have to dilate $P$ so that it
has an interior lattice point.  This is the \emph{co-degree} of the polytope.  \polymake
computes it with the command \texttt{LATTICE\_CODEGREE}.  In our case, this gives $1$, as
the cube contains the origin.  Recent results suggest that the degree is a more relevant
invariant of a lattice polytope than the dimension.  For instance, it is known that for
given degree $d$ and normalized volume $V$ or linear coefficient $h_1^*$ there is a
constant $c$ such that any lattice polytope of dimension $d\ge c$ is a lattice pyramid (a
pyramid where the apex is a lattice point with height $1$ over the
base)~\cite{BatyrevDegree, Nillhstar}.  The $h^*$-polynomials of a lattice polytope $P$
and the lattice pyramid with base $P$ coincide.

Finally, we can use \polymake to draw our cube, together with
its lattice points.  By default the command

\begin{figure}[tb]
  \centering
  \includegraphics[height=.15\textheight]{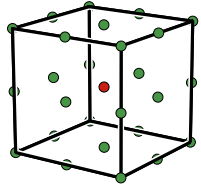}
  \caption{The  $3$-dimensional cube  with   its lattice  points.  The
    interior   lattice    points    are  drawn      in  a    different
    color.}
  \label{fig:cube}
\end{figure}

\begin{footnotesize}
\begin{verbatim}
polytope > $P->VISUAL->LATTICE_COLORED;
\end{verbatim}
\end{footnotesize}

triggers the visualization with \texttt{JavaView} \cite{javaview}.  See
Figure~\ref{fig:cube}.

\section{Analyzing an Example}

Let us look at the cone $C\subset\RR^6$ positively spanned by the rows of the $10\times
6$-matrix
\begin{equation}\label{eq:M}
  M \ = \
  \begin{pmatrix}
    0 & 1 & 0 & 0 & 0 & 0\\
    0 & 0 & 1 & 0 & 0 & 0\\
    0 & 0 & 0 & 1 & 0 & 0\\
    0 & 0 & 0 & 0 & 1 & 0\\
    0 & 0 & 0 & 0 & 0 & 1\\
    1 & 0 & 2 & 1 & 1 & 2\\
    1 & 2 & 0 & 2 & 1 & 1\\
    1 & 1 & 2 & 0 & 2 & 1\\
    1 & 1 & 1 & 2 & 0 & 2\\
    1 & 2 & 1 & 1 & 2 & 0
  \end{pmatrix} \, .
\end{equation}
The cone $C$ is \emph{pointed}, that is, it does not contain any line.  Equivalently, $C$
is projectively equivalent to a polytope $\bar C$.  The rows of $M$ are precisely the rays
(or \emph{generators}) of $C$, that is, they correspond to the vertices of $\bar C$.  The
key fact about $C$ is the following.

\begin{theorem}[Bruns et al.~\cite{BGHW99}]\label{thm:BGHMW}
  The vector $(9,13,13,13,13,13)$ lies in $C$, but it cannot be written as a non-negative
  integral linear combination of six generators of $C$.
\end{theorem}

This says that $C$ does not satisfy the \emph{integral Carath\'eodory property}, and thus
it is a counter-example to a conjecture of Seb\H{o} \cite{Sebo}.  We will sketch how this
can be verified using \polymake.  Moreover, we will reveal the combinatorial structure.
There is an integral transformation which maps $C$ to a cone with $0/1$-coordinates
\cite{BGHW99}, and there is also a realization of $\bar C$ as a lattice polytope.  Both
other representations could be used in the sequel with the same results.  The following
command creates a new matrix object representing the matrix $M$ above.  Here the user
types in the coefficient directly; alternatively, they could also be read from a file.

\begin{footnotesize}
\begin{verbatim}
polytope > $M=new Matrix<Rational>(<<".");
polytope (2)> 0 1 0 0 0 0
polytope (3)> 0 0 1 0 0 0
polytope (4)> 0 0 0 1 0 0
polytope (5)> 0 0 0 0 1 0
polytope (6)> 0 0 0 0 0 1
polytope (7)> 1 0 2 1 1 2
polytope (8)> 1 2 0 2 1 1
polytope (9)> 1 1 2 0 2 1
polytope (10)> 1 1 1 2 0 2
polytope (11)> 1 2 1 1 2 0
polytope (12)> .
\end{verbatim}
\end{footnotesize}

\polymake can work with pointed polyhedral cones right away, so it is
legal to write

\begin{footnotesize}
\begin{verbatim}
polytope > $C=new Polytope<Rational>(POINTS=>$M);
\end{verbatim}
\end{footnotesize}

$C$ in terms of the (rows of the) matrix $M$.  The first step is to verify that the
generators actually form the Hilbert basis of $C$.  Each integral cone admits a unique
minimal family of vectors such that any integral point inside can be written as a
non-negative linear combination of these.  Moreover, this family is finite, and this is
the \emph{Hilbert basis} of the cone.  \polymake cannot compute Hilbert bases directly,
but instead it relies of \normaliz \cite{normaliz2}, which uses an algorithm of Bruns and
Koch~\cite{BKalgo}.  The alternative implementation in \texttt{4ti2}~\cite{4ti2} uses a
lift and project approach described in~\cite{Hemmecke}.

\begin{footnotesize}
\begin{verbatim}
polytope > print $C->HILBERT_BASIS;
0 0 0 0 0 1
0 0 0 0 1 0
0 0 0 1 0 0
0 0 1 0 0 0
1 0 2 1 1 2
0 1 0 0 0 0
1 1 1 2 0 2
1 1 2 0 2 1
1 2 0 2 1 1
1 2 1 1 2 0
\end{verbatim}
\end{footnotesize} 

The output coincides with our first input, and this says that the generators of $C$ do
form a Hilbert basis.  One can show that it suffices to check if the vector
$x=(9,13,13,13,13,13)$ can be written as a non-negative integral linear combination of six
\emph{linearly independent} generators.  The following \polymake code enumerates all
possibilities.

\begin{footnotesize}
\begin{verbatim}
$x=new Vector<Rational>([9,13,13,13,13,13]);
foreach (all_subsets_of_k(6,0..9)) {
  $B=$M->minor($_,All);
  if (det($B)) {
    print lin_solve(transpose($B),$x), "\n";
  }
}
\end{verbatim}
\end{footnotesize} 

For each non-vanishing maximal minor $B$ we solve the linear system of equations $yB=x$,
and we print the unique solution to the screen.  The resulting 185 lines of output can be
checked by hand: All coefficients are integral, and each solution has at least one
negative coefficient.  Clearly, adding one or two more lines of code would also leave this
final check to \polymake.

In the remainder of this section we want to exploit \polymake's features to further
investigate the cone $C$ or rather the projectively equivalent polytope $\bar C$ from the
combinatorial point of view.  The first thing is to look at the facets (which had been
computed by \texttt{cddlib}~\cite{cddlib} before).  There are 27 of them.  Instead of
printing them all we only look at two, and instead of printing the coordinates we list the
numbers of the generators incident.

\begin{footnotesize}
\begin{verbatim}
polytope > print $C->VERTICES_IN_FACETS->[8];
{0 1 2 3 4}
polytope > print $C->VERTICES_IN_FACETS->[22];
{5 6 7 8 9}
\end{verbatim}
\end{footnotesize}

This shows that $\bar C$ has two disjoint facets of five vertices each. Since $\dim\bar
C=5$ each facet is a $4$-polytope, and this shows that both facets must be simplices.  The
numbers of the facets depend on the sequence of the output of \texttt{cddlib}, but the
numbers of the vertices correspond to the matrix $M$ as defined above.  \polymake uses the
first coordinate to homogenize.  By looking at \eqref{eq:M} we see that the first five
points have a leading zero coordinate, and hence the facet numbered~8 is the face at
infinity of $C$.  There is another popular $5$-polytope which happens to be the joint
convex hull of two disjoint $4$-dimensional simplices, and this is the $5$-dimensional
cross polytope.

\begin{footnotesize}
\begin{verbatim}
polytope > $cross5 = cross(5);
polytope > print isomorphic($C->GRAPH->ADJACENCY,$cross5->GRAPH->ADJACENCY);
1
\end{verbatim}
\end{footnotesize} 

The vertex-edge graph of $\bar C$ turns out to be isomorphic (as an abstract graph) to the
graph of the cross polytope.  This has been verified by \polymake's interface to \nauty
\cite{nauty}.  In fact, one can show that both polytopes even share the same $2$-skeleton.
If we compare the $f$-vectors we see that the cross polytope has five more facets and five
more \emph{ridges} (faces of codimension $2$) than $\bar C$.

\begin{footnotesize}
\begin{verbatim}
polytope > print $cross5->F_VECTOR - $C->F_VECTOR;
0 0 0 5 5
\end{verbatim}
\end{footnotesize}

This leads to a natural conjecture: What if, combinatorially, $\bar C$ can be constructed
from the cross polytope by picking five pairs of adjacent facets and ``straightening''
them?  Equivalently, the dual graph of $\bar C$ would result from the dual graph of the
cross polytope by contracting a partial matching of five edges.  This can be verified as
follows. First let us look at two more facets or rather the set of generators incident
with them.

\begin{footnotesize}
\begin{verbatim}
polytope > print $C->VERTICES_IN_FACETS->[12];
{0 2 5 7 8}
polytope > print $C->VERTICES_IN_FACETS->[13];
{1 2 5 7 8}
\end{verbatim}
\end{footnotesize}

The facets 12 and 13 with the vertices $\{0,2,5,7,8\}$ and $\{1,2,5,7,8\}$, respectively,
are adjacent in the dual graph (via the common ridge with vertex set $\{2,5,7,8\}$). This
edge and four others can be contracted in a copy of the dual graph of \verb+$cross5+.
Taking a copy first is necessary since \polymake's objects are immutable.

\begin{footnotesize}
\begin{verbatim}
polytope > $g=new props::Graph($cross5->DUAL_GRAPH->ADJACENCY);
polytope > $g->contract_edge(12,13);
polytope > $g->contract_edge(24,26);
polytope > $g->contract_edge(17,21);
polytope > $g->contract_edge(3,11);
polytope > $g->contract_edge(6,22);
polytope > $g->squeeze;
\end{verbatim}
\end{footnotesize}

It is only the last command which turns \verb+$g+ into a valid graph again.  The reason
for this is that \polymake's graphs necessarily have their nodes consecutively numbered.
Contracting an edge means to destroy one node (actually the second one).  Squeezing
renumbers the remaining vertices properly.

\begin{footnotesize}
\begin{verbatim}
polytope > print isomorphic($C->DUAL_GRAPH->ADJACENCY,$g);
1
\end{verbatim}
\end{footnotesize}

This final computation by \nauty explains the combinatorial structure of the cone $C$ in
Theorem~\ref{thm:BGHMW}, or the projectively equivalent polytope $\bar C$, completely.

\section*{Acknowledgements}
We are indebted to Ewgenij Gawrilow without whom \polymake would not exist.  Thanks to
Christian Haase and Benjamin Nill for valuable discussions. 

\bibliographystyle{amsplain} 
\providecommand{\bysame}{\leavevmode\hbox to3em{\hrulefill}\thinspace}
\providecommand{\MR}{\relax\ifhmode\unskip\space\fi MR }
\providecommand{\MRhref}[2]{%
  \href{http://www.ams.org/mathscinet-getitem?mr=#1}{#2}
}
\providecommand{\href}[2]{#2}

\end{document}